\documentclass[a4paper]{article}

\usepackage{graphics}
\usepackage{amsmath}
\usepackage{amssymb}
\usepackage{amscd}
\usepackage[ruled]{algorithm2e}

\usepackage[latin1]{inputenc}
\usepackage{varioref}
\usepackage{a4wide}

\newcommand{\rbn}{\ensuremath{R_b(n)}}
\newcommand{\rbnp}{\ensuremath{R_b(n+1)}}
\newcommand{\rbi}{\ensuremath{R_b(\infty)}}
\newcommand{\tbi}{\ensuremath{T_b(\infty)}}

\newcommand{\I}[1]{\ensuremath{I_{#1}}}
\newcommand{\C}[1]{\ensuremath{C_{#1}}}

\newcommand{\tr}[1]{\ensuremath{\stackrel{#1}{\longrightarrow}}}
\newcommand{\proof}{\textbf{Proof : }}
\newcommand{\cqfd}{\hfill \fbox{} \vskip 0.2cm}
\newcommand{\fb}[1]{\ensuremath{{\hookrightarrow_{#1}}}}
\newcommand{\inc}[2]{\ensuremath{#1^\fb{#2}}}

\newcommand{\Val}{v}

\newcommand{\fig}[2]{\begin{figure}[!h] \centerline{\includegraphics{#1.eps}}
\caption{\label{fig_#1} #2} \end{figure}}

\newtheorem{prop}{Proposition}
\newtheorem{theo}{Theorem}

\newtheorem{lemma}{Lemma}
\newtheorem{definition}{Definition}

\begin{document}

\begin{center}

{\bf \Large
Partitions of an Integer into Powers
}

\medskip

{\bf Matthieu Latapy}

\smallskip

\textsc{liafa}, Universit\'e Paris 7, 2 place Jussieu, 75005 Paris.

latapy@liafa.jussieu.fr

\end{center}


\begin{abstract}
In this paper, we use a simple discrete dynamical model
to study partitions of integers into powers of another integer.
We extend and generalize some known results about their enumeration
and counting, and we give new structural results.
In particular, we show that the set of these partitions can be ordered
in a natural way which gives the distributive lattice structure to this set.
We also give a tree structure which allow efficient and simple
enumeration of the partitions of an integer.
\end{abstract}

\section{Introduction}

We study here the problem of writing a non-negative integer $n$
as the sum of powers of another positive integer $b$:
$$n = p_0b^0 + p_1b^1 + \dots + p_{k-1}b^{k-1}$$
with $p_{k-1}\not= 0$ and $p_i \in \mathbb{N}$ for all $i$.
Following \cite{Rod69}, we call the $k$-tuple $(p_0,p_1,\dots,p_{k-1})$ 
a \emph{$b$-ary partition} of $n$.
The integers $p_i$ are called the \emph{parts} of the partition
and $k$ is the \emph{length} of the partition.
A $b$-ary partition of $n$ can be viewed as a representation
of $n$ in the basis $b$, with digits in $\mathbb{N}$. Conversely,
given a $k$-tuple $(p_0,\dots, p_{k-1})$ and a basis $b$, we will
denote by $\Val_b(p_0,\dots,p_{k-1})$ the integer
$p_0b^0 + p_1b^1 + \dots + p_{k-1}b^{k-1}$.
There is a unique $b$-ary partition such that $p_i < b$ for all $i$,
and it is the usual (canonical) representation of $n$
in the basis $b$. Here, we consider the problem without any
restriction over the parts: $p_i \in \mathbb{N}$, which is actually
equivalent to say that $p_i \in \lbrace 0, 1, \dots, n \rbrace$
for all $i$. We will mainly be concerned with the enumeration
and counting of the $b$-ary partitions of $n$, for given integers
$n$ and $b$.

This natural combinatorial problem has been introduced
by Mahler \cite{Mah40}, who showed that the logarithm of
the number of $b$-ary
partitions of $n$ grows as $\frac{(\log n)^2}{2\log b}$. This 
asymptotic approximation was later improved by de Bruijn \cite{Bru48} and
Pennington \cite{Pen53}. Knuth \cite{Knu66} studied the special
case where $b=2$.
In this case, the function counting the $b$-ary partitions for a
given $n$ is called the \emph{binary partition function}.
This function has been widely studied. Euler and Tanturri
\cite{Eul50,Tan18a,Tan18b}
studied its exact computation and Churchhouse \cite{Chu69,Chu71} studied
its congruence properties, while Fr\"oberg \cite{Fro77} gave
a final solution
to its asymptotical approximation. Later, R\"odseth \cite{Rod69}
generalized some of these results to $b$-ary partitions for any $b$.
Finally, Pfaltz \cite{Pfa95}
studied the subcase of the binary partitions of integers which
are powers of two.

We are concerned here with the exact computation of the
number of $b$-ary partitions of a given integer $n$, for any $b$.
We will use a powerful technique we developped in \cite{LP99}
and \cite{LMMP98}: incremental construction of the set of
$b$-ary partitions of $n$, infinite extension and coding by
an infinite tree. This method gives a deep understanding
of the structure of the set
of $b$-ary partitions of $n$. We will obtain this way a tree structure
which permits the enumeration of all the $b$-ary partitions of $n$ in
linear time with respect to their number. We will also
order these partitions in a natural way which gives
the distributive lattice structure to this set. We recall that 
a {\em lattice} is a partially ordered set such that
any two elements $a$ and $b$ have a least
upper bound (called {\em supremum} of $a$ and $b$ and denoted by $a \vee b$)
and a greatest lower bound (called {\em infimum} of $a$ and $b$ and denoted by
$a \wedge b$). The element $a \vee b$ is the smallest element among the
elements greater than both $a$ and $b$. The element $a \wedge b$ is
defined dually. A  lattice
is \emph{distributive} if for all $a$, $b$ and $c$:
$(a \vee   b)\wedge (a \vee   c) = a \vee   (b \wedge c)$ and 
$(a \wedge b)\vee   (a \wedge c) = a \wedge (b \vee   c)$.
A distributive lattice is a strongly structured set, and many
general results, for example efficient coding and algorithms,
are known about such sets. For more details,
see for example \cite{DP90}.

Notice that if we consider $b=1$ and restrict the problem to partitions of
length at most $n$, then we obtain the \emph{compositions} of $n$, i.e.
the series of at most $n$ integers, the sum of which equals $n$. Many studies
already deal with this special case. In particular, the (infinite)
distributive lattice $R_1(\infty)$ which we will introduce in
Section~\ref{sec_rbi} is isomorphic to the well known Young lattice
\cite{Ber71}. Therefore, we will suppose $b > 1$ in the following.
Notice however that some of the results we present here are
already known in
this special case (for example the distributive lattice structure),
therefore they can be seen as an extension of the existing ones.

\section{The lattice structure}
\label{sec_lat}

In this section, we define a simple dynamical model which generates
\emph{all} the $b$-ary partitions of an integer.
We will show that the set of $b$-ary
partitions, ordered by the reflexive and transitive closure of
the successor relation, has the distributive lattice structure.

Let us consider a $b$-ary partition $p = (p_0, p_1, \dots, p_{k-1})$
of $n$, and let
us define the following transition (or rewriting) rule:
$p \tr{i} q$ if and only if for all
$j \not\in \lbrace i,i+1 \rbrace$, $q_j = p_j$, $p_i \ge b$,
$q_i = p_i - b$ and $q_{i+1} = p_{i+1}+1$ (with the assumption
that $p_k = 0$). In other words, if $p_i$ is at least equal
to $b$ then $q$ is obtained from $p$ by removing $b$ units from
$p_i$ and adding one unit to $p_{i+1}$. We call this operation
\emph{firing $i$}. The important point is to notice
that $q$ is then a $b$-ary partition of $n$. We call $q$ a
\emph{successor}\,\footnote{Notice that the term \emph{successor} can have
many different meanings. We follow here the standard usage in discrete
dynamical models, but in order theory the term has another meaning,
and one may also consider that a \emph{successor} of an integer $n$ should
be the integer $n+1$, which is not the case here.}
of $p$, and we denote by $Succ_b(p)$ the set of
all the successors of $p$, with respect to the rule.
We denote by \rbn\ the set of $b$-ary partitions of $n$ reachable from $(n)$
by iterating the evolution rule, ordered by the reflexive and transitive
closure of the successor relation.
Notice that the successor relation is the covering
relation of the order, since it is defined as the transitive and
reflexive closure of the successor relation, and one can easily
verify that this relation has no reflexive
($x \longrightarrow x$) and no transitive ($x \longrightarrow z$ with
$x \longrightarrow y$ and $y \longrightarrow z$) edge.
See Figure~\ref{fig_exs_tr} for some examples.

\fig{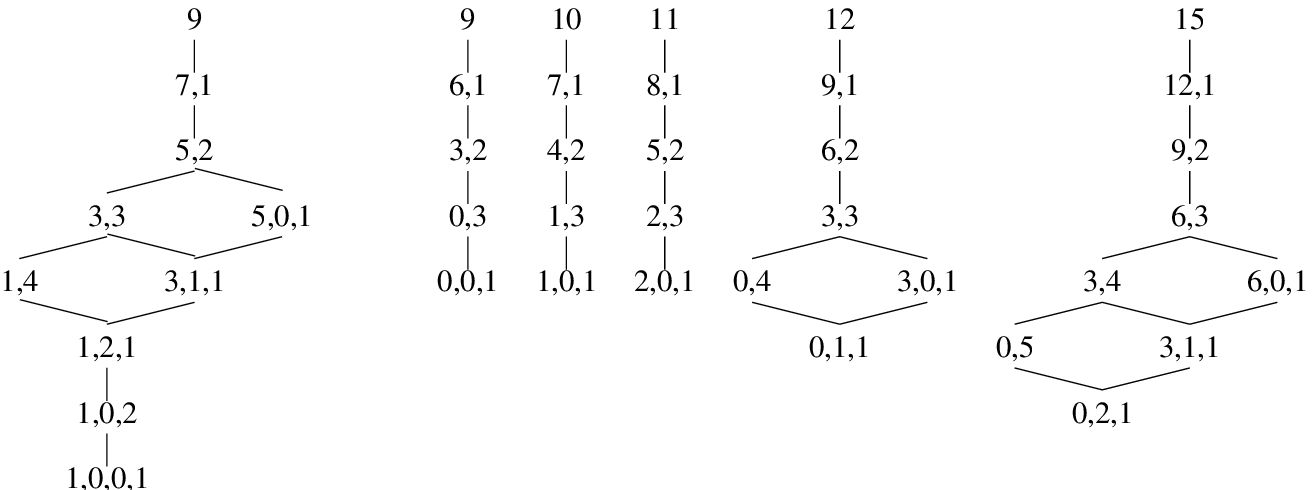}{From left to right, the sets $R_2(9)$, $R_3(9)$,
             $R_3(10)$, $R_3(11)$, $R_3(12)$ and $R_3(15)$.
             From Theorem~\ref{th_lat}, both of these sets is a
             distributive lattice.}


Given a sequence $f$ of firings, we denote by $|f|_i$ the number of
firings of $i$ during $f$.
Now, consider an element $p$ of \rbn, and two sequences
$f$ and $f'$ of firings which transform $(n)$ into $p$.
Then, $p_i = |f|_{i-1} - b \cdot |f|_i = |f'|_{i-1} - b \cdot |f'|_i$.
Suppose that there exists an integer $i$ such that $|f|_i \not= |f'|_i$,
and let $i$ be the smallest such integer. Then, $|f|_{i-1} = |f'|_{i-1}$
and the equality $|f|_{i-1} - b \cdot |f|_i = |f'|_{i-1} - b \cdot |f'|_i$
is impossible. Therefore, we have $|f|_i = |f'|_i$ for all $i$.
This leads to the definition
of the \emph{shot vector} $s(p)$: $s(p)_i$ is the number of times
one have to fire $i$ in order to obtain $p$ from $(n)$.
Now we can prove:

\begin{lemma}
\label{lem_ordre}
For all $p$ and $q$ in \rbn, $p \le q$ if and only if
for all $i$, $s(p)_i \ge s(q)_i$.
\end{lemma}
\proof
If $p \le q$, i.e. $p$ is reachable from $q$ then it is clear
that for all $i$, $s(p)_i \ge s(q)_i$.
Conversely, if there exists $i$ such that $s(p)_i > s(q)_i$,
then let $j$ be the smallest such integer. Therefore, $q_j > p_j + b$
and so $q$ can be fired at $j$. By iterating this process, we finally
obtain $p$, and so $p \le q$.
\cqfd

\begin{theo}
\label{th_lat}
For all integers $b$ and $n$, the order \rbn\ is
a \emph{distributive lattice} which contains \emph{all} the
$b$-ary partitions of $n$, with the infimum and supremum of
any two elements $p$ and $q$ defined by:
$$s(p\vee q)_i   = \min(s(p)_i,s(q)_i) \mbox{ and }
  s(p\wedge q)_i = \max(s(p)_i,s(q)_i).$$
\end{theo}
\proof
We first show that \rbn\ contains all the $b$-ary partitions of $n$.
Consider $p$ a $b$-ary partition of $n$. If $p=(n)$, then $p\in\rbn$,
so we suppose that $p \not= (n)$. Therefore, there must be an integer
$i>0$ such that $p_i > 0$. Let us define $q$ such that $q_j = p_j$
for all $j \not\in \lbrace i-1,i \rbrace$, $q_{i-1} = p_{i-1}+b$
and $q_i = p_i - 1$. It is clear that $q$ is a $b$-ary partition of
$n$, and that if $q \in \rbn$ then $p \in \rbn$ since
$q \tr{i-1} p$. It is also obvious that, if we iterate this
process, we go back to $(n)$, and so $p \in \rbn$.

We now prove the formula for the infimum and the supremum. Let
$p$ and $q$ be in $\rbn$, and $r$ such that $s(r)_i = \min(s(p)_i,s(q)_i)$.
From Lemma~\ref{lem_ordre}, $p$ and $q$ are reachable from $r$. Moreover,
if $p$ and $q$ are reachable from $t \in \rbn$, then, from
Lemma~\ref{lem_ordre}, $r$ is reachable
from $t$ since we must have $s(t)_i \le \min(s(p)_i,s(q)_i)$
(else one can not transform $t$ into $p$ or $q$). Therefore,
$r$ is the supremum of $p$ and $q$, as claimed in the theorem.
The argument for the infimum is symmetric. Finally, to prove
that the lattice is \emph{distributive}, we only have to check
that the formulae satisfy the distributivity laws.
\cqfd

We will now show that the dynamical model defined here can be viewed as a
special Chip Firing Game (CFG). A CFG \cite{BLS91,BL92} is defined over
a directed multigraph. A configuration of the game is a repartition of
a number of chips over the vertices of the graph, and it obeys the following
evolution
rule: if a vertex $\nu$ contains as many chips as its outgoing
degree $d$, then one can transfer one chip along each of its outgoing edges.
In other words, the number of chips at $\nu$ is decreased by $d$ and,
for each vertex $v \not= \nu$, the number of chips at $v$ is increased
by the number of edges from $\nu$
to $v$. This model is very general and has
been introduced in various contexts, such as physics, computer
science, economics, and others. It is in particular very close to the
famous Abelian Sandpile Model \cite{LP00}.

It is
known that the set of reachable configurations of such a game,
ordered with the reflexive and transitive closure of the transition
rule,
is a Lower Locally Distributive (LLD) lattice (see \cite{Mon90}
for a definition and properties), but it is not
distributive in general \cite{BL92,LP00,MPV01}. However, if a lattice
is LLD and its dual, i.e. the lattice obtained by reversing the order
relation, also is LLD, then the lattice is distributive. Therefore, we can
give another proof of the fact that \rbn\ is a distributive lattice by
showing that it is the set of reachable configurations of a CFG, and
that its dual too\,\footnote{This idea is due to Cl\'emence Magnien,
who introduced this new way to prove that a set is a distributive lattice
using two Chip Firing Games.}.

Given two integers $n$ and $b$, let us consider the
following multigraph $G=(V,E)$ defined by:
$V = \lbrace 0, \dots, n \rbrace$ and there are $b^{i+1}$ edges from
the $i$-th vertex to the $(i+1)$-th, for all $n < i \le 0$.
Now, let us consider the CFG $C$ defined over $G$ by the initial configuration
where the vertex $0$ contains $n$ chips, the other ones being empty.
Now, given a configuration $c$ of the CFG, where $c_i$ denotes the number of
chips in the vertex number $i$, let us denote
by $\bar{c}$ the vector such that $\bar{c}_i = \frac{c_i}{b^i}$.
Then, if the CFG is in the configuration $c$, an application of the rule
to the vertex number $i$ gives the configuration $c'$ such
that $c'_i = c_i - b^{i+1}$, $c'_{i+1} = c_{i+1} + b^{i+1}$ and
$c'_j = c_j$ for all $j \not\in \lbrace i, i+1 \rbrace$.
Notice that this means exactly that $\bar{c}_i$ is decreased by $b$ and
that $\bar{c}_{i+1}$ is increased by $1$, therefore
an application of the CFG rule corresponds exactly
to an application of the evolution rule we defined above, and so
the set of reachable configurations of the CFG is isomorphic to \rbn.
This leads to the fact that \rbn\ is a LLD lattice.

Conversely, let $G'$ be the multigraph obtained from $G$ by reversing each
edge, and let us consider the CFG $C'$ over $G'$ such that the initial
configuration of $C'$ is the final configuration of $C$. Then it is clear
that the set of reachable configurations of $C'$ is nothing
but the dual of the one of $C$, therefore it is isomorphic to the
dual of \rbn. This leads to the fact that the dual of \rbn\ is a
LLD lattice, which allows us to conclude that \rbn\ is a distributive
lattice.

\section{From \rbn\ to \rbnp}

In this section, we give a method to construct the transitive reduction
(i.e. the successor relation) of \rbnp\ from
the one of \rbn. In the following, we will simply call this the
\emph{construction of \rbnp\ from \rbn}.
This will show the self-similarity of these sets, and give
a new way, purely structural, to obtain a recursive
formula for $|\rbn|$, which is previously known from \cite{Rod69}
(the special case where $b=2$ is due to Euler \cite{Eul50}).
This construction will also show the special role played by certain
$b$-ary partitions, which will be widely used in the rest of the
paper. Therefore, we introduce a few notations
about them. We denote by $P_i(b,n)$ the set of the partitions
$p$ in \rbn\ such that $p_0 = p_1 = \dots = p_{i-1} = b-1$. Notice
that for all $i$ we have $P_i(b,n) \subseteq P_{i+1}(b,n)$
and that $P_0(b,n)=R_b(n)$.
If $p=(p_0,\dots,p_{k-1})$ is in $P_i(b,n)$, we denote
by \inc{p}{i}\ the $k$-uple $(0,\dots,0,p_i+1,p_{i+1},\dots,p_{k-1})$.
In other words, \inc{p}{i} is obtained from $p$ by switching all
the $i$ first components of $p$ from $b-1$ to $0$ and adding one unit to
its $i$-th componend\,\footnote{This operator is known in numeration studies
as an odometer. See~\cite{GLT95} for more precisions.}.
Notice that the $k$-uple $\inc{p}{0}$,
which is simply obtained from $p$ by adding one unit to its
first component, is always a $b$-ary partition of $n+1$.
If $S$ is a subset of $P_i(b,n)$,
we denote by \inc{S}{i}\ the set $\lbrace \inc{p}{i} |\ p\in S\rbrace$.

Notice that, if $p \tr{i} q$ in \rbn, then $\inc{p}{0} \tr{i} \inc{q}{0}$
in \rbnp. This remark makes it possible to construct
\rbnp\ from \rbn: the construction procedure starts with
the lattice \inc{\rbn}{0}\ given by its diagram. Then, we look for those
elements in \inc{\rbn}{0}\ that have a successor out of \inc{\rbn}{0}.
The set of these elements will be denoted by \I{0}, with
$\I{0} \subseteq \inc{\rbn}{0}$.
At this point, we add all the missing successors of the elements of $I_0$.
The set of these new elements will be denoted by \C{0}. Now, we
look for the elements in \C{0} that have a successor out of the
constructed set. The set of these elements is denoted by \I{1}.
More generally, at the $i$-th step of the procedure we
look for the elements in \C{i-1} with missing
successors and call \I{i} the set of these elements. We add the new
successors of the elements of \I{i} and call the set of these new elements
\C{i}. At each step, when we add a new element, we also add its covering
relations. Since \rbnp\ is a finite set, this procedure terminates.
At the end, we obtain the whole set \rbnp. In the rest of this section,
we study more precisely this construction process.

\begin{lemma}
\label{lem_succ}
Let $p$ be a $b$-ary partition in $P_i(b,n)$. If $p_i\not= b-1$
then $Succ_b(\inc{p}{i}) = \inc{{Succ_b(p)}}{i}$. Else,
$Succ_b(\inc{p}{i}) = \inc{{Succ_b(p)}}{i} \cup \lbrace \inc{p}{i+1} \rbrace$.
\end{lemma}
\proof
If a transition $p \tr{j} q$ is possible, then
$\inc{p}{i} \tr{j} \inc{q}{i}$ is obviously possible. Moreover, an
additional transition is possible from \inc{p}{i}\ if and only if
$p_i=b-1$. In this case, $\inc{p}{i} \tr{i} \inc{p}{i+1}$.
\cqfd

\begin{lemma}
\label{lem_Pi}
\label{lem_bij}
For all integer $b$, $n$ and $i$, we define the function
$r_i: P_i(b,n) \rightarrow R_b(\frac{n+1}{b^i}-1)$
by: $r_i(p)$ is obtained from $p \in P_i(b,n)$ by
removing its $i$ first components (which are equal to $b-1$).
Then, $r_i$ is a bijection.
\end{lemma}
\proof
Let us consider $p$ in $P_i(b,n)$:
$p = (b-1,b-1,\dots,b-1,p_i,\dots,p_k)$.
Then, it is clear that $r_i(p)=(p_i,\dots,p_k)$ is in
$R_b(\frac{n-(b-1)-(b-1)b-\dots-(b-1)b^{i-1}}{b^i}) =
R_b(\frac{n+1-b^i}{b^i}) = R_b(\frac{n+1}{b^{i+1}}-1)$.
Conversely, if we consider $p$ in $R_b(\frac{n+1}{b^i}-1)$,
then $r_i^{-1}(p)=(b-1,b-1,\dots,b-1,p_0,p_1,\dots,p_k)$ is
a $b$-ary partition
of $m = (b-1) + (b-1)b + \dots + (b-1)b^{i-1} + \frac{n+1-b^i}{b^i}$,
which is nothing but $n$. Therefore, $r_i^{-1}(p)$ is in \rbn.
\cqfd

\begin{lemma}
\label{lem_IC}
For all integer $b$, $n$ and $i$, we have 
$\I{i} = \inc{P_{i+1}(b,n)}{i}$ and $\C{i} = \inc{P_{i+1}(b,n)}{{i+1}}$.
\end{lemma}
\proof
By induction over $i$. For $i=0$, it is clear from
Lemma~\ref{lem_succ} that the set of elements
in \inc{\rbn}{0}\ with a missing successor, namely $I_0$, is exactly
\inc{P_1(b,n)}{0}. Moreover, the set of these missing successors, namely
$C_0$, is clearly \inc{P_1(b,n)}{1}. Now, let us suppose that the claim
is proved for $i$ and let us prove it for $i+1$. The set $I_{i+1}$ is
the set of elements in $C_i$ with one missing successor. By induction
hypothesis, we have $C_i = \inc{P_{i+1}(b,n)}{{i+1}}$ and so, from
Lemma~\ref{lem_succ}, $I_{i+1} = \inc{P_{i+2}(b,n)}{{i+1}}$. Then, by
application of the evolution rule, it is clear that the set $C_{i+1}$
of the missing successor is \inc{P_{i+2}(b,n)}{{i+2}}, which proves the
claim.
\cqfd

\begin{theo}
For any positive integer $b$ and $n$, we have:
$$R_b(n) = \bigsqcup_{i\ge0} r_i^{-1}\inc{\left(R_b\left(\frac{n}{b^i}-1\right)\right)}{i}$$
$$|R_b(n)| = \sum_{i=0}^{\lfloor n/b \rfloor} \begin{tabular}{|c|} $R_b(\frac{i}{b})$ \end{tabular}$$
where $\bigsqcup$ denotes the disjoint union, where \rbn\ is taken
as $\emptyset$ when $n$ is not a positive integer, and with
$R_b(0)=\lbrace 0 \rbrace$.
\end{theo}
\proof
From the construction procedure described above,
we have $\rbn = \inc{R_b(n-1)}{0} \sqcup \bigsqcup_{i\ge 0} C_i$.
From Lemma~\ref{lem_IC}, we obtain
$\rbn = \inc{R_b(n-1)}{0} \sqcup \bigsqcup_{i\ge 0} \inc{P_{i+1}(b,n)}{{i+1}}$.
Moreover, since $\inc{R_b(n-1)}{0}$ is nothing but $\inc{P_0(b,n)}{0}$,
this is equivalent to $\rbn = \bigsqcup_{i\ge 0} \inc{P_i(b,n)}{i}$.
Finally, from Lemma~\ref{lem_Pi}, we obtain the announced formula.

From this formula, we have $R_b(\frac{n}{b}) = 
 \bigsqcup_{i\ge0} r^{-1}(\inc{R_b(\frac{n}{b^{i+1}}-1)}{i})$.
 Therefore,
  $|\rbn| = \sum_{i\ge0} | R_b(\frac{n}{b^i}-1) |
          = |R_b(n-1)| + \sum_{i\ge0} | R_b(\frac{n}{b^{i+1}}-1) |
          = |R_b(n-1)| + | R_b(\frac{n}{b}) |$.
 We obtain the claim by iterating this last formula.
\cqfd

\begin{figure}
\begin{center}
\rotatebox{90}{
\begin{minipage}{20cm}
\resizebox{20cm}{12cm}{\includegraphics{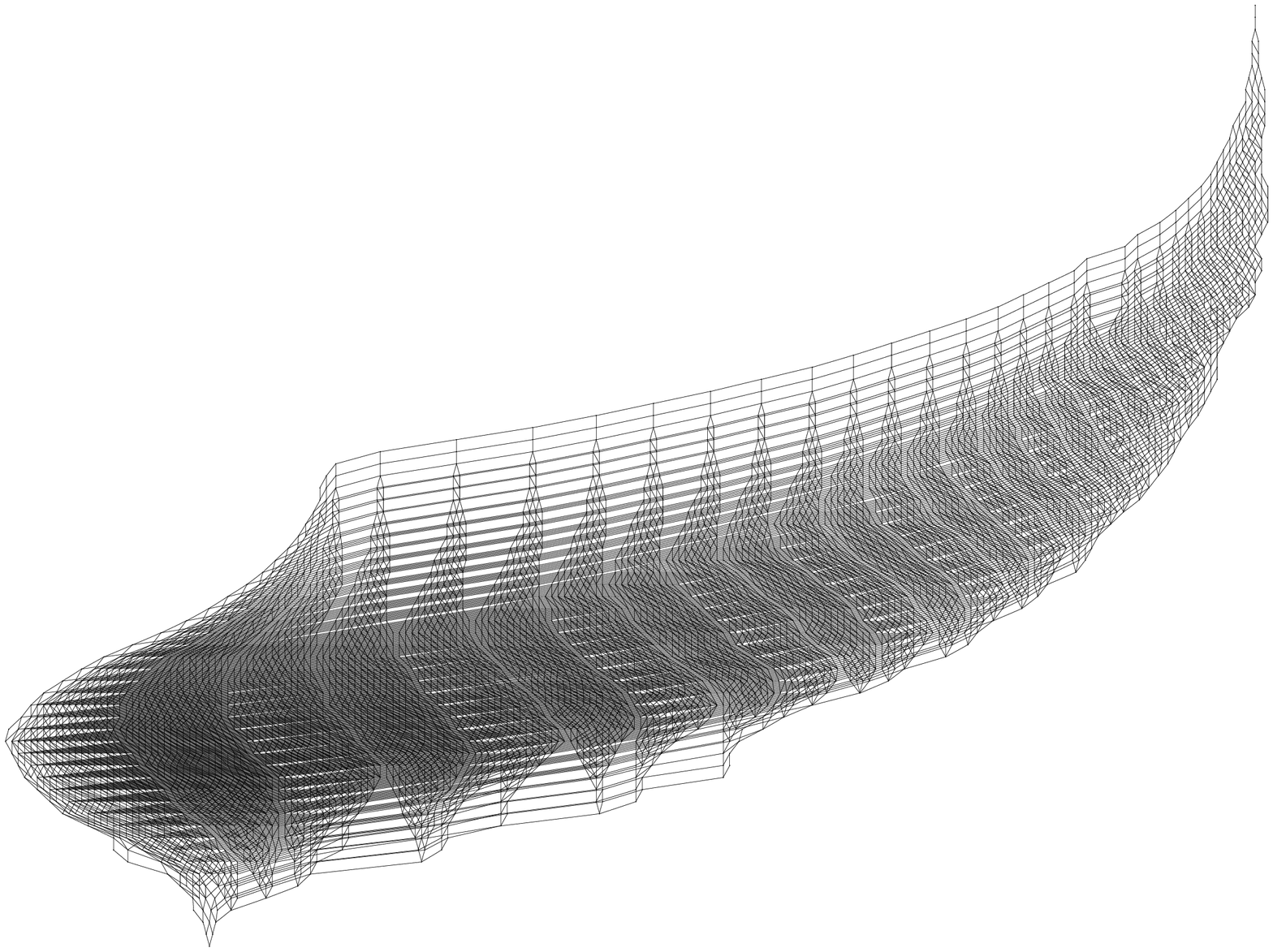}}
\caption{\label{fig_rb2_80}The distributive lattice $R_2(80)$, which
contains $4124$ elements and $12484$ edges.
The self-similarity of the set clearly appears on this diagram.}
\end{minipage}}
\end{center}
\end{figure}

The first formula given in this theorem can be used to compute the sets
\rbn\ efficiently since it only involves \emph{disjoint} unions. We will
give in Section~\ref{sec_tree} another method to compute \rbn\ which
is much simplier, as it gives \rbn\ a tree structure. However, the
formula is interesting since it points out the self-similar structure
of the set (see Figure~\ref{fig_rb2_80}).

The second formula is previouly
known from \cite{Rod69}, and from \cite{Eul50} in the special
case where $b=2$. Notice that this does not give a way to compute
$|\rbn|$ in linear time with respect to $n$, which is an unsolved
problem in the general case, but it gives a very simple way to
compute recursively $|\rbn|$.

\section{Infinite extension}
\label{sec_rbi}

\rbn\ is the lattice of the $b$-ary partitions of $n$ reachable from
$(n)$ by iteration of the evolution rule. We
now define \rbi\ as the set of all $b$-ary partitions reachable
from $(\infty)$. The order on \rbi\ is
the reflexive and transitive closure of the successor relation.
For $b=2$, the first $b$-ary
partitions in \rbi\ are given in Figure~\ref{fig_rb_infi} along with
their covering relation (the first component, which is always
infinity, is not represented on this diagram).
Notice that it is still possible to define the shot vector $s(p)$
of an element $p$ of \rbi\ by: $s(p)_i$ is the number of times one has
to fire $i$ in order to obtain $p$ from $(\infty)$.

\fig{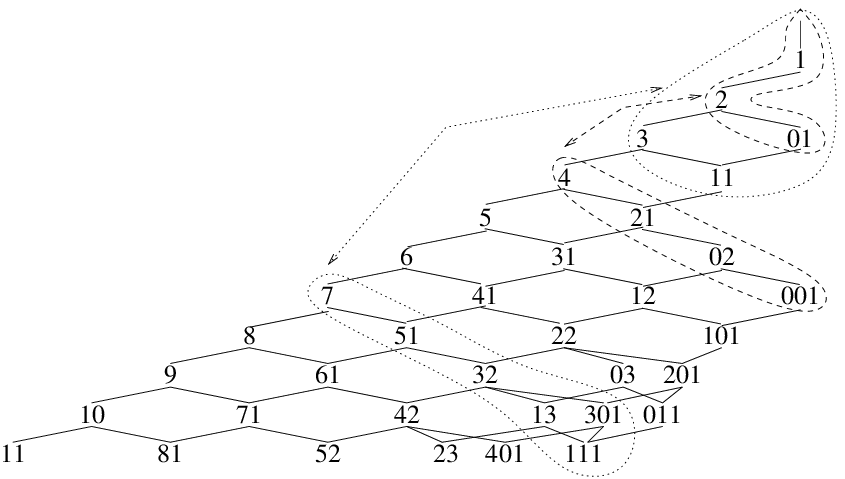}{The first $b$-ary partitions obtained in \rbi\ when $b=2$.
Two parts isomorphic to $R_2(4)$ are distinguished, as well as two parts
isomorphic to $R_2(7)$.}

\begin{theo}
The set \rbi\ is a distributive lattice with:
$$s(p\vee q)_i   = \min(s(p)_i,s(q)_i) \mbox{ and }
  s(p\wedge q)_i = \max(s(p)_i,s(q)_i)$$
for all $p$ and $q$ in \rbi.
Moreover, for all $n$ the functions
$$
\pi : s=(s_1, s_2,\dotsi,s_k) \longrightarrow
\pi(s) = (\infty,s_2,\dots,s_k)
$$
and
$$
\tau : s=(s_1, s_2,\dotsi,s_k) \longrightarrow
\tau(s) = (\infty,s_1,s_2,\dots,s_k)
$$
are lattice embeddings of \rbn\ into \rbi.
\end{theo}
\proof
The proof for the distributive lattice structure and for the formulae
of the infimum and supremum is very similar to the proof of
Theorem~\ref{th_lat}. Therefore, it is left to the reader.

Given $p$ and $q$ in \rbn, we now prove that
$\pi(p) \vee \pi(q) = \pi(p \vee q)$. From Theorem~\ref{th_lat}, we have
$s(p \vee q)_i = \min(s(p)_i,s(q)_i)$. Moreover, it is clear that
$s(\pi(x))_i = s(x)_i$ for all $x$ in \rbn. Therefore,
$s(\pi(p \vee q))_i = \min(s(\pi(p))_i,s(\pi(q))_i))$, which
shows that $\pi$ preserves the supremum. The proof of
$\pi(p) \wedge \pi(q) = \pi(p \wedge q)$ is symmetric. Therefore,
$\pi$ is a lattice embedding.

The proof for $\tau$ is very similar
when one has noticed that the shot vector of $\tau(s)$ is obtained
from the one of $s$ by adding a new first component equal to $n$.
\cqfd

With similar arguments, one can easily show that $\pi(\rbn)$ is
a sublattice of $\pi(\rbnp)$, and so we have an infinite chain of
distributive lattices:
$$\pi(R_b(0)) \le \pi(R_b(1)) \le \dots \le \pi(R_b(n)) \le\pi(R_b(n+1)) \le \dots \le \rbi,$$
where $\le$ denotes the sublattice relation.
Moreover, one can use the self-similarity estalished here to construct
filters of \rbi\ (a \emph{filter} of a poset is an upper closed part of
the poset). Indeed, if one defines $R_b(\le n)$ as the sub-order of
\rbi\ over $\cup_{i\le n}R_b(i)$, then one can construct efficiently
$R_b(\le n+1)$ from $R_b(\le n)$ by extracting from $R_b(\le n)$ a part
isomorphic to \rbnp\ and pasting it to $R_b(\le n)$.
See Figures~\ref{fig_rb_infi} and \ref{fig_rb2_80}.


Notice that, for all integer $b$, \rbi\ contains exactly
all the finite sequences of integers, since any such sequence can be viewed
as a $b$-ary partition of an integer $n$. Therefore, we provide
infinitely many ways to give the set of finite sequences of integers the
distributive lattice structure.

\section{Infinite tree}
\label{sec_tbi}
\label{sec_tree}

As shown in our construction of \rbnp\ from \rbn, each $b$-ary partition
$p$ in \rbnp\ is obtained from another one $p' \in \rbn$ by
application of the \inc{}{}\ operator: $p = \inc{\mbox{$p'$}}{i}$
with $i$ an integer between $0$ and $l(p')$, where $l(p')$ denotes
the number of $b-1$ at the beginning of $p'$. Thus, we can define an
infinite tree \tbi\ whose nodes are
the elements of $\bigsqcup_{n\ge 0}{\rbn}$ and in which
the fatherhood relation is defined by:
$$
q \mbox{ is the $(i+1)$-th son of $p$ if and only if } q = \inc{p}{i}
\mbox{ for some } i,\ 0 \le i \le l(p).
$$
The root of this tree is $(0)$ and each node $p$ of \tbi\ has $l(p)+1$ sons.
The first levels of \tbi\ when $b=2$ are shown in Figure~\ref{fig_tb_infi}
(we call the set of elements of depth $n$ the ``level $n$''
of the tree).

\fig{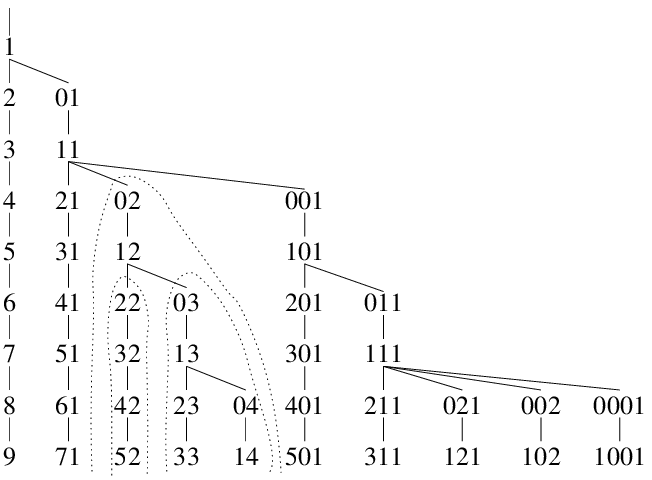}{The first levels of \tbi\ when $b=2$. We distinguished some
special subtrees, which will play an important role in the following.}

\begin{prop}
The level $n$ of \tbi\ contains exactly the elements of \rbn.
\end{prop}
\proof
Straightforward from the construction of \rbnp\ from \rbn\ given above
and the definition of the tree.
\cqfd

\noindent
If we define $\overline{\rbn}$ as $\lbrace (s_2, \dots, s_k) \ |\ 
(s_1,s_2,\dots,s_k) \in \rbn\rbrace$, then:

\begin{prop}
For all integer $n$, the elements of $\overline{\rbn}$ are exactly
the elements of the $\lfloor \frac{n}{b} \rfloor$ first levels of \tbi.
\end{prop}
\proof
Let us first prove that the elements of \rbn\ are the nodes of a subtree
of \tbi\ that contains its root. This is 
obviously true for
\mbox{$n=0$}. The general case follows by induction, since
by construction
the elements of $\overline{\rbnp} \setminus \overline{\rbn}$ are
sons of elements of $\overline{\rbn}$.

Now, let us consider an element $e$ of the $l$-th level of \tbi. If there is a 
$b$-ary partition $p$ of $n$ such that $\overline{p} = e$, then
clearly $p_i = e_{i-1}$ for all $i > 0$ and $p_0 = n - b \cdot l$.
Therefore, if $e$ is in $\overline{\rbn}$ then all the elements of
the $l$-th level are in $\overline{\rbn}$, and this is clearly
the case exactly when $0 \le l < \lfloor \frac{n}{b} \rfloor$.
This ends the proof.
\cqfd

Notice that this proposition gives a simple way to enumerate the
elements of \rbn\ for any $n$ in linear time with respect to
their number, since it gives this set a tree structure.
Algorithm~\ref{algo_enum} acheives this.

\begin{algorithm}
\KwIn{An integer $n$ and a basis $b$}
\KwOut{The elements of \rbn}
\Begin{   
 $\mbox{Resu} \leftarrow \lbrace (n) \rbrace$\;
 $\mbox{CurrentLevel} \leftarrow \leftarrow \lbrace () \rbrace$\;
 $\mbox{OldLevel} \leftarrow \emptyset$;\ 
 $l \leftarrow 0$\;
 \While{$l < \lfloor \frac{n}{b} \rfloor$}{
  $\mbox{OldLevel} \leftarrow \mbox{CurrentLevel}$\;
  $\mbox{CurrentLevel} \leftarrow \emptyset$\;
  $l \leftarrow l + 1$\;
  \ForEach{$p$ in OldLevel}{
   $i \leftarrow 0$\;
   \Repeat{$p_{i-1} \not= b-1$}{
    Add $\inc{p}{i}$ to CurrentLevel\;
    $i \leftarrow i+1$\;
   }
  }
 \ForEach{$e$ in CurrentLevel}{
  Create $p$ such that $p_i = e_{i-1}$ for all $i>0$ and $p_0 = n - b\cdot l$\;
  Add $p$ to Resu\;
  }
 }
 Return(Resu)\;
}
\caption{\label{algo_enum}Efficient enumeration of the elements of \rbn.}
\end{algorithm}

We will now show that \tbi\ can be described recursively, which allows us
to give a new recursive formula for $|\rbn|$. In order to do this, we
will use a series known as the $b$-ary carry sequence \cite{Slo73}:
$c_b(n) = k$ if $b^k$ divides $n$ but $b^{k+1}$ does not.
Notice that this function is defined only for $n>0$ (or one can consider
that $c_b(0) = \infty$).
These series appear in many contexts, and have many equivalent
definitions\,\footnote{
For example, if one defines the series $C_{b,0} = 0$
and $C_{b,i} = C_{b,i-1},
\stackrel{\mbox{$b-1$ times}}{\overbrace{i,C_{b,i-1}}}$,
then $c_b(i)$ is nothing but the $i$-th integer of the series
$C_{b,i}$. The ten first values for $c_2(i)$ are $0,1,0,2,0,1,0,3,0,1$
and the ten first ones for $c_3(i)$ are $0,0,1,0,0,1,0,0,2,0$.}.
Here, we will mainly use the fact that the
first $n$ such that $c_b(n)=k$ is $n=b^k$, and the fact that $c_b(n)$ is
nothing but the number of components equal to $b-1$ at the begining of
the canonical representation of $n-1$ in the basis $b$.

\begin{definition}
Let $p \in \tbi$. Let us consider the rightmost branch of \tbi\ rooted
at $p$ ($p$ is considered as the first node of the branch).
We say that $p$ is the root of a \emph{$X_{b,k}$ subtree (of \tbi)} if
this rightmost branch is as follows:
for $i \le b^{k-1}$,
the $i$-th node on the branch has $j=c_b(i) + 1$ sons,
and the $l$-th ($1 \le l < j$) of these sons is the root
of a $X_{b,l}$ subtree.
Moreover, the $(b^{k-1}+1)$-th node of the branch is itself
the root of a $X_{b,k}$ subtree.
\end{definition}

For example, we show in Figure~\ref{fig_tb_infi} a $X_{2,2}$ subtree of
$T_2(\infty)$, composed of a $X_{2,1}$ subtree and another $X_{2,2}$ subtree.
Notice that a $X_{b,1}$ subtree is simply a chain.

\begin{prop}
\label{prop_Xk}
Let $p=(0,0,\dots,0,p_k,\dots)$ in \tbi\ with $p_k > b-1$.
Then, $p$ is the root of a $X_{b,k+1}$ subtree of \tbi.
\end{prop}
\proof
The proof is by induction over $k$ and the depth of $p$.
Let us consider the rightmost branch
rooted at $p$. Since, for all $q$ in \tbi, the rightmost son
of $q$ is \inc{q}{i}\ with $i$ the number of $b-1$ at the beginning
of $q$, it is clear that the $j$-th node of this branch for $j \le b^k$ is
$q=(q_0,\dots,q_{k-1},p_k,\dots)$ where $(q_0,\dots,q_{k-1})$ is
the canonical representation of $j-1$ in the basis $b$.
Therefore, $q$ begins with $c_b(j)$ components equal to $b-1$, and so,
for $l=1,\dots,c_b(j)$, the $l$-th son of $q$ starts with
$l-1$ zeroes followed by a component equal to $b > b-1$.
By induction hypothesis, we then have that the sons
of $q$ are the roots of $X_{b,l}$ subtrees.
Moreover, the $(b^k+1)$-th node on the rightmost branch
begins with exactly $k$ zeroes followed by a component greater than $b-1$,
and so it is the root of a $X_{b,k+1}$ subtree by induction hypothesis.
\cqfd

\begin{theo}
\label{th_tbi}
The infinite tree \tbi\ is a $X_{b,\infty}$ tree: it is
a chain (its rightmost branch) such that its $i$-th node
has $c_b(i)$ sons and the $j$-th of these sons, $1 \le j \le c_b(i)$, is the
root of a $X_{b,j}$ subtree. Moreover, the $i$-th node of
the chain is the canonical representation of $i-1$ in
the basis $b$.
\end{theo}
\proof
Since the rightmost son of $p \in \tbi$ is \inc{p}{i},
where $i$ is the number of $b-1$ at the beginning of $p$,
and since the root of \tbi\ is nothing but the canonical
representation of $0$,
it is clear by induction that the $i$-th node of the rightmost branch of
\tbi\ is the canonical representation of $i-1$ in the basis $b$.
Then, the theorem follows from Proposition~\ref{prop_Xk}.
\cqfd

We now have a recursive description of \tbi, which allows us to
give recursive formula for the cardinal of some special sets.
Let us denote by $\pi_b(l,k)$ the number of paths of length exactly $l$
starting from the root of a $X_{b,k}$ subtree of \tbi. We have:

\begin{theo}
$$
\pi_b(l,k) = \left\{\begin{array}{ll}
1 & \mbox{if $0 \le l < b$}\\
1 + \sum_{i=1}^{l} \sum_{j=1}^{c_b(i)} \pi_b(l-i,j) & \mbox{if $b \le l \le b^{k-1}$}\\
\pi_b(l-b^{k-1},k) + \sum_{i=1}^{b^{k-1}} \sum_{j=1}^{c_b(i)} \pi_b(l-i,j) & \mbox{otherwise ($l > b^{k-1}$)}
\end{array}
\right.
$$
Moreover, $|\rbn| = \pi_b(n,n)$ and the number of $b$-ary partitions of
$n$ into exactly $l$ parts is $\pi_b(n-(b-1)^l,l)$.
\end{theo}
\proof
The formula for $\pi_b(l,k)$ is directly deduced from the definition
of the $X_{b,k}$ subtrees. The other formulae derive from
Theorem~\ref{th_tbi} and from the fact that all the $b$-ary partitions
of length $l$ are in a $X_{b,l}$ subtree of \tbi which is rooted at the
$(b-1)^l$-th node of the righmost branch of \tbi.
\cqfd

\section{Perspectives}

The results presented in this paper mainly point out
the strong self-similarity and the structure of the
sets \rbn. As already noticed, it is an open question
to compute the cardinal of \rbn\ in linear time with
respect to $n$, and one may expect to obtain a solution
using these results.

Another interesting direction is to investigate how one
can extend the dynamics we study. A first idea is to consider
non-integer basis, in particular complex basis or Fibonnacci basis.
For example, if we consider the complex basis $b = i-1$
then we can obtain all the
ways to write an integer $n$ as the sum of powers of $b$ by iterating
the following evolution rule from $(n)$: $q$ is a successor of $p$
if $p - q = (0, \dots, 0,2,0,-1,-1,0\dots,0)$.
In other words, we can decrease by two the $j$-th component
of $p$ and increase by one its $(j+2)$-th and its $(j+3)$-th
components for some integer $j$. This gives to the set of representations
of $n$ in the complex basis $b=i-1$ the lattice structure,
since this can be encoded by a Chip Firing Game \cite{LP00}
(notice however that in this case the lattice is no longer
distributive). Another interesting case is when $b=1$. As already
noticed, we obtain the Young lattice, or equivalently the lattice
of the compositions of $n$.


\section{Acknowledgments}
I thank Christiane Frougny and Cl\'emence Magnien
for many useful comments on preliminary versions,
which deeply improved the manuscript.

\bibliographystyle{alpha}
\bibliography{bib2021,biblio}

\end{document}